\documentclass[journal=iecred, manuscript=article, table]{achemso}

\usepackage{chemformula} 
\usepackage[T1]{fontenc} 
\usepackage{amsfonts}
\usepackage{amsmath}
\usepackage{amssymb}
\usepackage{dsfont}
\usepackage{mathtools}
\usepackage{upgreek}
\usepackage{textcomp} 
\usepackage{siunitx}
\usepackage{graphicx}
\usepackage[version = 4]{mhchem}
\usepackage{algorithm}
\usepackage{algpseudocode}
\usepackage{array}
\usepackage{booktabs}
\usepackage{multicol}
\usepackage{multirow}
\usepackage[flushleft]{threeparttable}
\usepackage{placeins}  
\usepackage[final]{changes}
\usepackage{etoolbox}
\usepackage{framed}
\usepackage{nomencl} 

\setkeys{acs}{maxauthors=99}

\def\ve#1{\mathchoice{\mbox{\boldmath$\displaystyle#1$}}
	{\mbox{\boldmath$\textstyle#1$}}
	{\mbox{\boldmath$\scriptstyle#1$}}
	{\mbox{\boldmath$\scriptscriptstyle#1$}}} 

\algnewcommand{\Inputs}[1]{%
  \State \textbf{Inputs:}
  \Statex \hspace*{\algorithmicindent}\parbox[t]{.8\linewidth}{\raggedright #1}
}
\algnewcommand{\Initialize}[1]{%
  \State \textbf{Initialize:}
  \Statex \hspace*{\algorithmicindent}\parbox[t]{.8\linewidth}{\raggedright #1}
}
\algnewcommand{\ReturnList}[1]{%
  \State \textbf{Return:}
  \Statex \hspace*{\algorithmicindent}\parbox[t]{.8\linewidth}{\raggedright #1}
}

\makenomenclature
\renewcommand\nomgroup[1]{%
  \item[\bfseries
  \ifstrequal{#1}{M}{Main Symbols}{%
  \ifstrequal{#1}{S}{Subscripts and Supercripts}{%
  \ifstrequal{#1}{A}{Acronyms}{}%
  }}%
]}

\makeatletter
\@ifundefined{chapter}
  {\def\wilh@nomsection{section}}
  {\def\wilh@nomsection{chapter}}

\def\thenomenclature{%
  \begin{multicols}{2}[%
    \csname\wilh@nomsection\endcsname*{\nomname}
    \if@intoc\addcontentsline{toc}{\wilh@nomsection}{\nomname}\fi
    \nompreamble]
  \list{}{%
    \labelwidth\nom@tempdim
    \leftmargin\labelwidth
    \advance\leftmargin\labelsep
    \itemsep\nomitemsep
    \let\makelabel\nomlabel}%
}
\def\endthenomenclature{%
  \endlist
  \end{multicols}
  \nompostamble}
\makeatother

\newcolumntype{x}[1]{>{\centering\arraybackslash\hspace{0pt}}m{#1}}

\DeclareMathOperator\erf{erf}

\definechangesauthor[name={Benoit Chachuat}, color=orange]{BC}
\definechangesauthor[name={Kennedy P. K.}, color=blue]{KPK}
\setaddedmarkup{\uwave{#1}}

\author{Kennedy P. Kusumo}
\affiliation{Centre for Process Systems Engineering, Department of Chemical Engineering, Imperial College London, London, SW7 2AZ, UK}
\author{Lucian Gomoescu}
\affiliation{Centre for Process Systems Engineering, Department of Chemical Engineering, Imperial College London, London, SW7 2AZ, UK}
\alsoaffiliation{Process Systems Enterprise Ltd, London, W6 7HA, UK}
\author{Radoslav Paulen}
\affiliation{Faculty of Chemical and Food Technology, Slovak University of Technology in Bratislava, 812 43 Bratislava, Slovakia}
\author{Salvador Garc\'ia Mu\~noz}
\affiliation{Small Molecule Design and Development, Lilly Research Laboratories, Eli Lilly \& Company, Indianapolis, IN 46285, USA}
\author{Constantinos~C.~Pantelides}
\affiliation{Centre for Process Systems Engineering, Department of Chemical Engineering, Imperial College London, London, SW7 2AZ, UK}
\alsoaffiliation{Process Systems Enterprise Ltd, London, W6 7HA, UK}
\author{Nilay Shah}
\affiliation{Centre for Process Systems Engineering, Department of Chemical Engineering, Imperial College London, London, SW7 2AZ, UK}
\author{Beno\^it Chachuat}
\affiliation{Centre for Process Systems Engineering, Department of Chemical Engineering, Imperial College London, London, SW7 2AZ, UK}
\email{b.chachuat@imperial.ac.uk}

\title[Bayesian Approach to Probabilistic Design Space Characterization]{Bayesian Approach to\\ Probabilistic Design Space Characterization:\\ A Nested Sampling Strategy}


\begin{document}

\begin{tocentry}
\includegraphics[width=\linewidth]{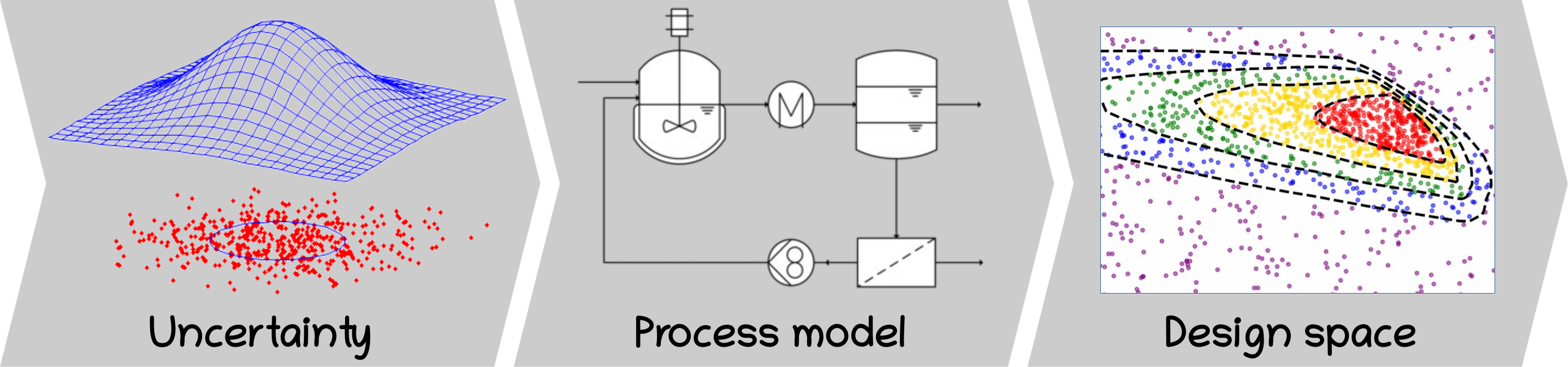}
\end{tocentry}

\begin{abstract}
Quality by design in pharmaceutical manufacturing hinges on computational methods and tools that are capable of accurate quantitative prediction of the design space. This paper investigates Bayesian approaches to design space characterization, which determine a feasibility probability that can be used as a measure of reliability and risk by the practitioner. An adaptation of nested sampling---a Monte Carlo technique introduced to compute Bayesian evidence---is presented. The nested sampling algorithm maintains a given set of live points through regions with increasing probability feasibility until reaching a desired reliability level. It furthermore leverages efficient strategies from Bayesian statistics for generating replacement proposals during the search. Features and advantages of this algorithm are demonstrated by means of a simple numerical example and two industrial case studies. It is shown that nested sampling can outperform conventional Monte Carlo sampling and be competitive with flexibility-based optimization techniques in low-dimensional design space problems. Practical aspects of exploiting the sampled design space to reconstruct a feasibility probability map using machine learning techniques are also discussed and illustrated. Finally, the effectiveness of nested sampling is demonstrated on a higher-dimensional problem, in the presence of a complex dynamic model and significant model uncertainty.
\end{abstract}

\section{Introduction}

Over the years the pharmaceutical industry has identified an increasing need for systematic and holistic approaches to drug development and manufacturing, which has led to a high penetration of process systems engineering (PSE) tools\cite{troup2013,reklaitis2017}. This need is prompted by safety concerns and regulations alongside growing pressure to increase efficiency, both in production and in process development. A recent estimate~\cite{dimasi2016innovation} amounts to US\$2.6Bn for a single novel drug put on the market. To improve practices in the industry the International Council for Harmonization of Technical Requirements for Pharmaceuticals for Human Use (ICH) introduced the Quality by Design (QbD) initiative through a series of five guidelines: ICH Q8--12\cite{ICH}. Historically, pharmaceutical development and manufacturing had emphasized checklist-based operations rather than scientific understanding. So, by promoting a scientific and risk-based approach to pharmaceutical product development and manufacturing, the QbD initiative and ICH guidelines triggered a paradigm shift in the industry and a complete new range of activities for the practitioners.

\nomenclature[A]{PSE}{Process Systems Engineering}
\nomenclature[A]{QbD}{Quality by Design}
\nomenclature[A]{ICH}{International Council for Harmonization}

The QbD approach defines a Quality Target Product Profile (QTPP), which is a prospective summary of the quality characteristics of the pharmaceutical product that ensures the desired quality, safety and efficacy. The (physical, chemical, biological or microbiological) properties that should be within an appropriate limit, range or distribution to ensure the desired product quality are called Critical Quality Attributes (CQAs). \deleted{The limits of acceptability of such CQAs is known as the Design Space (DS). More precisely, } \replaced{A}{a} \replaced{Design Space (DS)}{DS} consists of {\em ``the multidimensional combination and interaction of input variables (material attributes) and process parameters that have been demonstrated to provide assurance of quality''}\cite{ICHQ8R2}. In practice a regulatory process of post-approval change is not required so long as the process parameters vary within the limits of the approved DS.

\nomenclature[A]{QTPP}{Quality Target Product Profile}
\nomenclature[A]{DS}{Design Space}
\nomenclature[A]{CQA}{Critical Quality Attribute}

A classical approach to DS characterization entails the following four steps\cite{boukouvala2010design,chatterjee2017}: (i) conduct a thorough experimental investigation of the relationships between process parameters and CQAs; (ii) analyze the sensitivity of the process parameters on the CQAs to only retain those parameters presenting a significant sensitivity; (iii) establish a mathematical or graphical representation of the DS using data-driven modeling and optimization techniques; and (iv) test the final DS by running further validation experiments, prior to submitting it for approval to the regulatory agencies. Multivariate statistical techniques, e.g. based on latent-variable modeling, have been proposed to reduce the cost and time needed to conduct these experiments\cite{facco2015,bano2018probabilistic}. Nevertheless, the presence of a large number of external disturbances and/or potential control manipulations lead to high-dimensional DS, whose effective characterization can be empowered by the use of process models\cite{garcia2015definition}. Such process models may be either data-driven or knowledge-driven insofar as they describe the relationships between process parameters and CQAs accurately.

Characterizing a DS based on a process model is akin to the popular concepts of resilience and flexibility introduced by the PSE community during the 1980s\cite{biegler1997,grossmann2014}. Two classical problems in this literature are: (i) the flexibility test\cite{haleman1983}, which verifies that a feasible operation can be obtained for specific process parameters under a range of uncertainty scenarios---possibly by manipulating certain controls; and (ii) the flexibility index\cite{swaney1985indexpart1}, which represents the maximal deviation of the process parameters from given nominal values that can be tolerated in order for the operation to remain feasible under a range of uncertainty scenarios. A related concept is that of process operability\cite{georgakis2003,lima2010}, which starts from an estimated range of the process uncertainty and calculates the desired ranges of the control variables so that a control strategy can guarantee a feasible operation.

Like the result of any model-based computation, the reliability of any DS determined on the basis of a process model depends on the accuracy of that model. If the uncertain parameters of the model follow a certain probability distribution then the DS itself is not a well-delineated region but a probabilistic one, where each point in the process parameter space is characterized by a probability of meeting \replaced{all the established limits on CQAs}{the QTPP}\cite{adjiman2009designspace}. Interestingly, a closely related concept to that of probabilistic DS was introduced in the PSE literature later during the 1990s under the name stochastic flexibility\cite{straub1990,pistikopoulos1990}, as a way of measuring the probability of a process design to operate feasibly in the presence of uncertainty. This probabilistic interpretation is in agreement with \citeauthor{peterson2008bayesian}\cite{peterson2008bayesian} who defines the DS as the set of process parameters such that the manufactured product satisfies \replaced{all of the CQA limits}{the QTPP} alongside other process constraints with a probability greater than a given threshold.

Existing computational approaches to probabilistic DS characterization differ in how they account for process model uncertainty and how they approximate the DS itself. Process model uncertainty is represented in either one of two ways: (i) set-membership, most commonly a joint confidence ellipsoid at a given confidence level as derived from a frequentist parameter estimation;\cite{garcia2015definition,laky2019optimization} or (ii) sampled distribution, for instance the sampled parameter posterior in Bayesian estimation or a bootstrap parameter distribution\cite{peterson2008bayesian,peterson2017,bano2018probabilistic}. Caution must be exerted in both cases since assuming an ellipsoidal region for non-normally distributed parameters or representing a probability distribution with an insufficient number of samples can lead to large errors in the predicted DS. 

In terms of DS representation, one may broadly categorize the existing approaches as: (i) design-centering algorithms; and (ii) sampling algorithms. The former consider a class of parameterized shapes---most commonly a box---and transform the DS characterization into an optimization problem that seeks to inscribe the largest possible shape within the DS. The flexibility-index problem\cite{swaney1985indexpart1,straub1993} falls into this category and it was recently applied to probabilistic DS computation\cite{ochoa2018,laky2019optimization}. A drawback with this approach is that a simple shape may not provide a good approximation of the DS and thus introduce significant conservatism. Another drawback is that design-centering problems give rise to complex mathematical programs with either robust (semi-infinite) or chance constraints that are computationally hard to tackle rigorously\cite{floudas2001,harwood2017}. By contrast, sampling algorithms discretize the process parameter range and return a subset of the samples that satisfy \replaced{all of the CQA limits}{the QTPP} and other constraints up to the desired reliability value. Exhaustive sampling may be achieved via a fine uniform gridding or (quasi) Monte Carlo sampling. The assessment of each parameter sample can be made via the solution of an optimization problem, which is akin to a flexibility test and is especially suited to a set-membership description of the model uncertainty\cite{adjiman2009designspace,boukouvala2017,laky2019optimization}. If the model parameter uncertainty is represented by a sampled distribution instead, both Monte Carlo and Bayesian techniques can be used for propagating the uncertainty to the CQAs and estimate a feasibility probability\cite{peterson2017,bano2018probabilistic}. These sampling techniques have proven effective in practice but they are computationally expensive and mainly tractable for low-dimensional DS at present. For computationally demanding process models surrogate-based approaches can be applied, including kriging, radial basis functions, and high-dimensional model representation\cite{wang2018}. For instance, surrogate-based adaptive sampling\cite{boukouvala2012feasibility,rogers2015a,rogers2015b} has been shown to speed-up the computation of flexibility indices. But while it can be computationally cheaper to characterize the DS using a simple surrogate, one needs to account for the additional burden of constructing the surrogate itself as well as any additional approximation error that may be introduced. 

Despite the large body of research on DS computation, there is a clear need for algorithms with improved computational efficiency to tackle industrially-relevant problems that have more design parameters or greater uncertainty. The focus of this paper is on nested sampling, a Monte Carlo technique that was introduced by \citet{skilling2004nested} for estimating the Bayesian evidence in parameter estimation. This algorithm proceeds by progressively sampling in nested contours of increasing likelihood, so as to maintain a dense enough sample in regions of higher likelihood, in the manner of adaptive sampling. An advantage of nested sampling over Markov chain Monte Carlo (MCMC) techniques is its better ability to handle multimodal posteriors\cite{Feroz2009}.

\nomenclature[A]{MCMC}{Markov Chain Monte Carlo}

Our main contribution herein is an adaptation of nested sampling for the characterization of a probabilistic DS. In this context nested sampling determines a set of samples with a feasibility probability larger than a given reliability threshold. A by-product of the algorithm is a second set of samples with a feasibility probability below the desired threshold, which provides an approximation of the entire probabilistic DS. In this manner nested sampling can offer a greater flexibility in selecting the probability threshold to provide assurance of quality. An important motivation behind the adaptation of nested sampling is the availability of efficient strategies for generating replacement proposals as the algorithm progresses\cite{mukherjee2006nested, Feroz2009}, which we can leverage to the benefit of DS characterization. Another appeal is the ability to readily exploit a sampled joint posterior distribution of the model parameters: by combining the proposed approach with a Bayesian estimation procedure for the model parameters, one can arrive at a truly Bayesian approach to DS characterization directly from experimental data. Last but not least importantly nested sampling is, like other sampling-based techniques, non-intrusive in the sense that it relies on the result of model simulations at given process parameter values only. In principle, this allows for black-box models such as a process flowsheet or \added{a} CFD model to be used for DS characterization.

\nomenclature[A]{CFD}{Computational Fluid Dynamics}

The rest of the paper is organized as follows: In the following section we review the mathematical formulation of the design space characterization problem. Then we present the nested sampling approach and use a simple example to illustrate its main features. We also test the method on two case studies of increasing complexity and discuss the results, before concluding the paper.

\section{Problem Statement}

Consider a manufacturing process for a pharmaceutical product that has its quality defined by some CQAs, denoted by $\ve s\in\mathds R^{n_s}$. Assume that a mathematical model of the process (either knowledge- or data-driven) is available that predicts the CQAs corresponding to a set of process parameters, denoted by $\ve d\in\mathcal K$ within the \textsl{knowledge space} $\mathcal K\subset\mathds R^{n_d}$:
\begin{align}
\ve s = \ve f\left(\ve d, \ve \theta\right)
\label{eq:processmodel}
\end{align}
The model parameters, $\ve \theta \in \mathds{R}^{n_\theta}$ may represent physical constants, coefficients in a regression model, or disturbances that affect the CQAs. The mapping $\ve f:\mathds R^{n_d}\times\mathds R^{n_\theta}\to\mathds R^{n_s}$ needs not be given in closed-form but could be implicitly defined by a set of algebraic or differential equations or even a black-box functions such as a process simulator. Assume furthermore that the \replaced{CQA limits}{constraints on the CQAs} are represented alongside other process constraints by the following inequalities:
\begin{align}
\ve G\left(\ve d, \ve \theta\right) := \ve g\left(\ve d, \ve f\left( \ve d, \ve \theta \right) \right) \leq \ve 0
\label{eq:constraints}
\end{align}
Notice that the variable $\ve s$ can be abstracted away from the inequality constraints without loss of generality. Similarly to the mapping $\ve f$ above, the constraint function $\ve G:\mathds R^{n_d}\times\mathds R^{n_\theta}\to\mathds R^{n_g}$ is not necessarily available in closed form. 

\nomenclature[M]{$s$}{critical quality attribute}
\nomenclature[M]{$d$}{process parameter}
\nomenclature[M]{$\theta$}{model parameter}
\nomenclature[M]{$\mathcal K$}{knowledge space}
\nomenclature[S]{nom}{nominal}
\nomenclature[M]{$\mathcal D$}{design space}
\nomenclature[M]{$G$}{critical quality attribute constraint}
\nomenclature[M]{$\mathds{P}[\cdot]$}{probability}
\nomenclature[M]{$\alpha$}{reliability value}
\nomenclature[M]{$p(\cdot)$}{posterior distribution}
\nomenclature[M]{$\mathcal{N}(\mu, \sigma)$}{Normal distribution with mean $\mu$ and standard deviation $\sigma$}

Ignoring the uncertainty in the model parameters leads to a nominal design space:
\begin{align}
\label{eq:nom_ds_definition}
\mathcal D_{\rm nom} := \{ \ve d\in \mathcal K: \ve G (\ve d, \ve \theta_{\rm nom}) \leq \ve 0 \}
\end{align}
for a given nominal value $\ve \theta_{\rm nom}$ of the model parameters. However, the value of $\ve \theta$ is inherently uncertain in practice by the nature of the modeling exercise. A Bayesian framework considers $\ve\theta$ as random variables with a joint distribution $p(\ve\theta)$ that describes the belief on the value of $\ve\theta$. For instance, $p(\ve\theta)$ could be estimated from experimental data using Bayesian inference. In this framework the model may only be used to predict the probability that the manufacturing process is feasible for a given $\ve d\in \mathcal{K}$, which is akin to a stochastic flexibility test\cite{straub1990}:
\begin{align}
\label{eq:ful_prob_def}
\mathds{P} \left[\ve G(\ve d,\cdot) \leq \ve 0\: \middle|\: p(\ve\theta)\right] :=
\int_{\{\ve \theta : \ve G(\ve d,\ve \theta) \leq \ve 0\}} p(\ve \theta)\, {\rm d}\ve \theta
\end{align} 
The problem of interest throughout this paper is to determine the probabilistic DS given by:
\begin{align}
\label{eq:prob_ds_definition}
\mathcal D_\alpha := \left\{ \ve d\in \mathcal K: \mathds{P}\left[ \ve G(\ve d,\cdot) \leq \ve 0 \: \middle|\: p(\ve\theta) \right] \geq \alpha \right\}
\end{align}
where $0<\alpha\leq 1$ is the so-called reliability value \cite{peterson2008bayesian}. Notice that the following set-membership counterpart to the probabilistic DS is often computed instead of $\mathcal D_\alpha$ in practice\cite{laky2019optimization}:
\begin{align}
\label{eq:robust_ds_definition}
\widehat{\mathcal D}_\alpha := \left\{ \ve d\in \mathcal K: \forall \ve\theta\in\Theta_\alpha, G(\ve d,\ve \theta) \leq \ve 0 \right\}
\end{align}
with $\Theta_\alpha$ chosen as the highest posterior density (HPD) set such that $\int_{\ve\theta\in\Theta_\alpha} p(\ve\theta)=\alpha$. However, $\mathcal D_\alpha$ and $\widehat{\mathcal D}_\alpha$ are not equivalent in general when the constraints are nonlinear or the model parameters are not normally distributed. \deleted{It is also worth noting an important assumption underlying the definition of $\mathcal D_\alpha$ (or $\widehat{\mathcal D}_\alpha$), which is that the process model is either structurally correct or that any structural mismatch can be reported in terms of parametric uncertainty}

\nomenclature[A]{HPD}{Highest Posterior Density}

According to the monotonicity property that $\mathcal D_{\alpha} \supset \mathcal D_{\alpha'}$ whenever $\alpha<\alpha'$, a higher reliability value increases conservatism by shrinking the DS. A practical choice for $\alpha$ entails trading-off the risk of violating \replaced{the CQA limits}{QTPP} against the loss of operational flexibility: a lower $\alpha$ value increases false positives in $\mathcal{D}_\alpha$ by including design parameters that will not fulfill all of the \replaced{CQA limits}{QTPP} and process constraints as expected; while a higher $\alpha$ value increases false negatives by excluding design parameters that are in fact feasible. Because false positives pose a threat to the assurance of quality, practitioners are prompted to be conservative. But the choice of $\alpha$ remains specific to the nature of the process at hand and to the risks involved. 

\added{It is also worth noting an important assumption underlying the definition of $\mathcal D_\alpha$ (or $\widehat{\mathcal D}_\alpha$), namely that the process model is either structurally correct or that any structural mismatch can be reported in terms of parametric uncertainty. A practitioner's lacking confidence in the process model structure can always increase conservatism by opting for a higher reliability value $\alpha$. However, the probabilities conveyed by $\mathcal{D}_\alpha$ could become misleading in the presence of large structural mismatch. A possible remedy entails the consideration of multiple candidate models under the Bayesian framework and the restatement of $\mathds{P} \left[\ve G(\ve d,\cdot) \leq \ve 0\: \middle|\: p(\ve\theta)\right]$ to explicitly consider structural uncertainty.}


\subsection{Illustrative Example}

We consider a simple case study with two design variables, $\ve d:=(d_1, d_2)^\intercal$ and a single CQA, $s$. We assume the following relationship between the design variables and the CQA:
\begin{align}
\label{eq:model_ds_example}
s =\ & \theta d_1^2 + d_2
\end{align}
with the model parameter $\theta$. The goal is to characterize the probabilistic DS inside the knowledge space $\mathcal{K} := [-1, 1]^2$ imposed by the following \replaced{CQA limits}{QTPP}:
\begin{equation}
\label{eq:CQA_ds_example}
0.20 \leq s \leq 0.75
\end{equation}

We consider the nominal case to be given by $\theta_{\rm nom}=1$, which gives the following nominal design space:
\begin{align}
\label{eq:nom_ds_example}
\mathcal D_{\rm nom} := \{ \ve d\in [-1,1]^2: 0.20\leq d_1^2+d_2 \leq 0.75 \}
\end{align}
Under a normality assumption for the model parameters, $\theta \sim \mathcal N(\mu_\theta,\sigma_\theta)$ the probabilistic DS can be expressed analytically as:
\begin{align}
\label{eq:prob_ds_example}
\mathcal D_\alpha :=\ & \left\{ \ve d\in [-1,1]: \mathds{P}\left[0.20\leq \theta d_1^2+d_2 \leq 0.75 \:\middle|\: \theta \sim \mathcal N(\mu_\theta,\sigma_\theta) \right] \geq \alpha \right\}\\
=\ & \left\{ \ve d\in [-1,1]: \erf\left(\frac{0.75-\mu_\theta d_1^2-d_2}{\sqrt{2}\sigma_\theta d_1^2}\right) - \erf\left(\frac{0.20-\mu_\theta d_1^2-d_2}{\sqrt{2}\sigma_\theta d_1^2}\right) \geq 2\alpha \right\} \nonumber
\end{align}

\begin{figure}[tb]
\centering
\includegraphics[width=.48\textwidth]{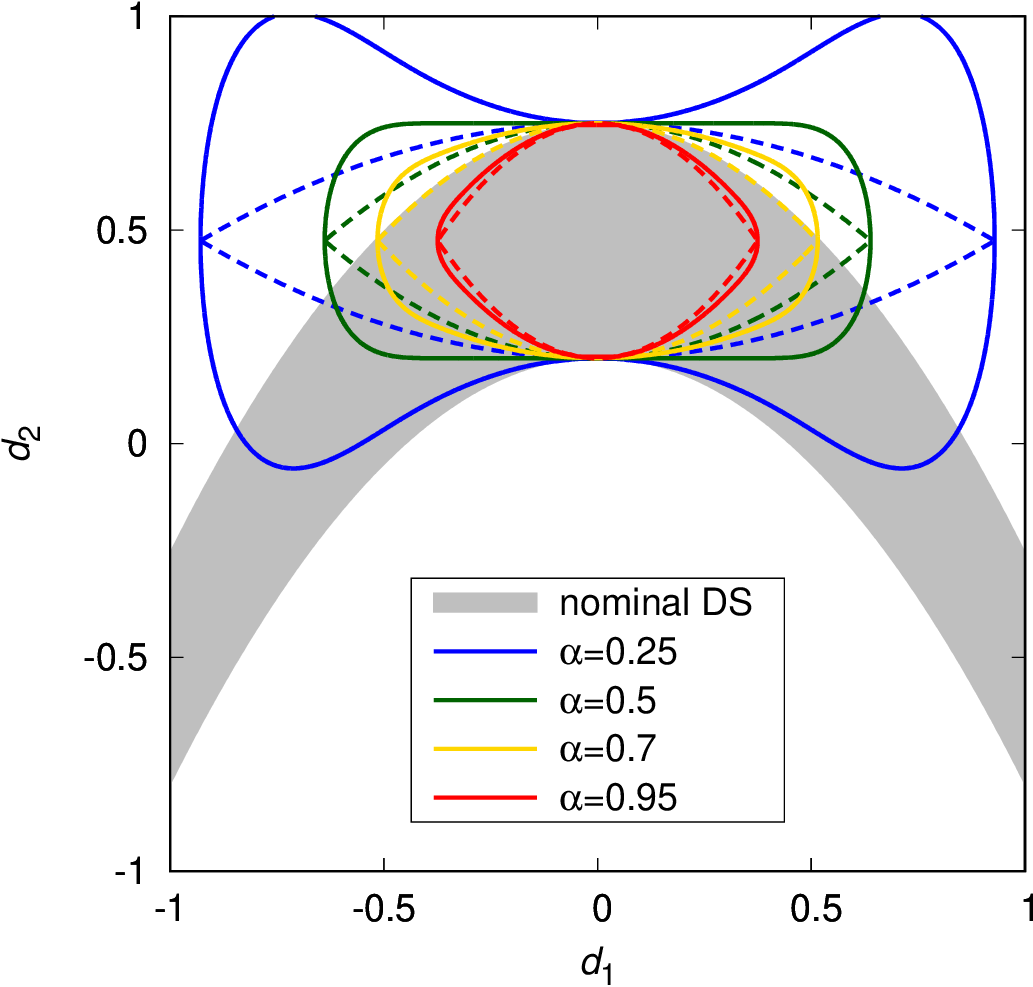}\hfill
\includegraphics[width=.48\textwidth]{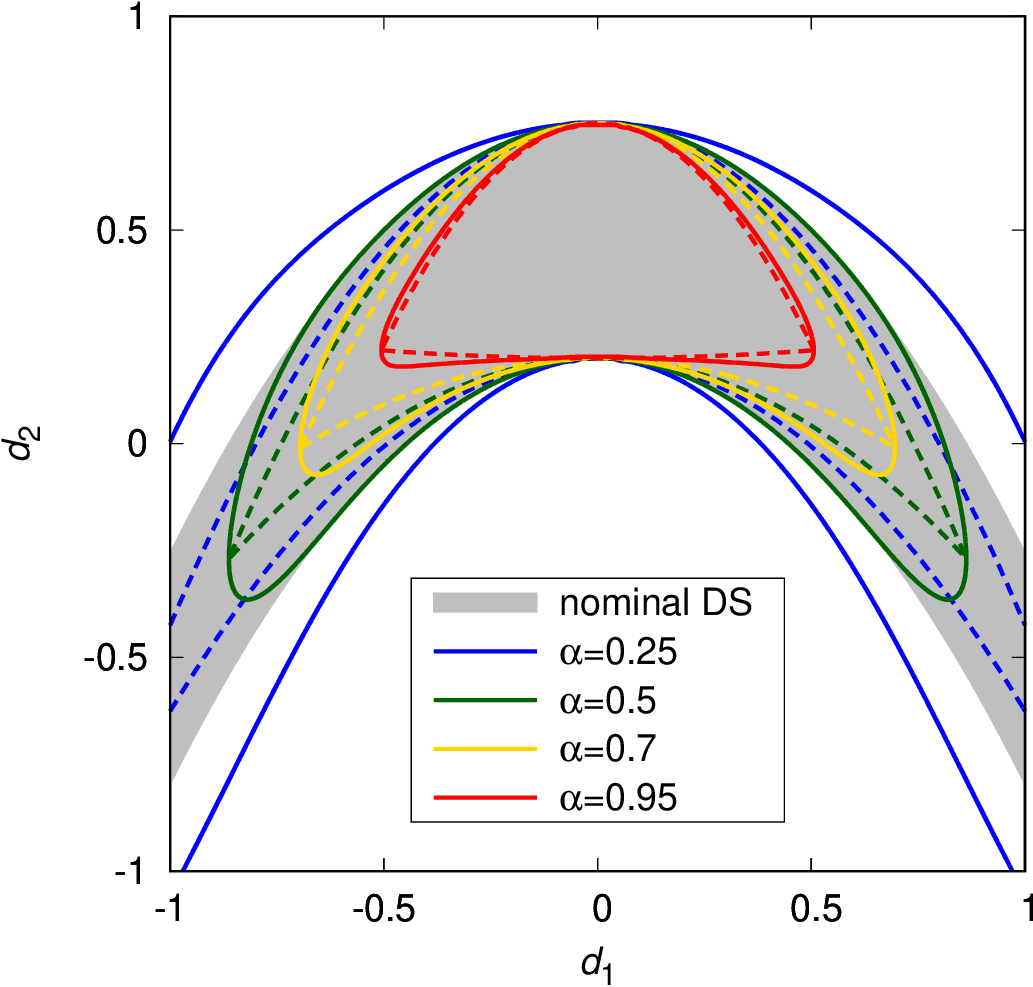}
\caption{Depiction of the design space for the illustrative example. The nominal design space in Equation~\eqref{eq:nom_ds_example} is represented with the gray-shaded area. The probabilistic design spaces corresponding to different reliability values $\alpha$ in Equation~\eqref{eq:prob_ds_example} are shown in solid lines with different colors. The set-membership counterparts in Equation~\eqref{eq:robust_ds_example} are shown in dashed lines. The left and right plots are for the uncertainty scenarios $\theta \sim \mathcal N(0,1)$ and $\theta \sim \mathcal N(1,\sqrt{0.3})$, respectively.}
\label{fig:comp_DS_example}
\end{figure}

Figure~\ref{fig:comp_DS_example} compares $\mathcal D_\alpha$ for two parameter uncertainty scenarios, $\theta \sim \mathcal N(0,1)$ in the left plot and $\theta \sim \mathcal N(1,\sqrt{0.3})$ in the right plot, and different reliability values $\alpha$. As expected, a more precise and accurate description of the model uncertainty enables a more favorable trade-off \replaced{between CQA}{QTPP} satisfaction and operational flexibility. With an imprecise and inaccurate uncertainty description (left plot), avoiding false positives requires a high reliability value (e.g. $\alpha=0.95$), but this results in a small DS and many false negatives. As the precision and accuracy of the model parameter improves (right plot), a lower reliability value may be selected without generating false positives (e.g. $\alpha=0.5$) and at the same time enable a much larger DS. 

For comparison, Figure~\ref{fig:comp_DS_example} also shows the DS counterpart for a set-membership uncertainty description based on Equation~\eqref{eq:robust_ds_definition}, which can also be expressed analytically as:
\begin{align}
\label{eq:robust_ds_example}
\widehat{\mathcal D}_\alpha :=\ & \left\{ \ve d\in [-1,1]: \forall \theta\in\mu_\theta\pm z_\alpha\sigma_\theta, 0.20\leq \theta d_1^2+d_2 \leq 0.75 \right\}\\
=\ & \left\{ \ve d\in [-1,1]: 0.20-(\mu_\theta-z_\alpha\sigma_\theta) d_1^2 \leq d_2 \leq 0.75-(\mu_\theta+z_\alpha\sigma_\theta) d_1^2 \right\} \nonumber
\end{align}
with the quantile $z_\alpha=\sqrt{2}\erf^{-1}(\alpha)$. For the same reliability value $\alpha$, $\widehat{\mathcal D}_\alpha$ is generally more conservative than $\mathcal D_\alpha$. And the mismatch between $\mathcal D_\alpha$ and $\widehat{\mathcal D}_\alpha$ furthermore reduces as $\alpha$ gets closer to $1$. The rest of the paper focuses on computational approach to characterizing the Bayesian design space $\mathcal D_\alpha$.

\section{Methodology}

Our main focus is on computational methods that rely on sampling to characterize a probabilistic design space. These methods are appealing in that they allow for the use of a process model in explicit or implicit form and they can accommodate any probability distribution for the model parameters. A pseudo-code for the standard Monte Carlo approach is presented in Algorithm~\ref{alg:MC}, where both the process parameter space range and the model uncertainty distribution are discretized. At each process parameter sample, $\ve d_i\in\mathcal S_d$ the probability of \replaced{CQA}{QTTP} satisfaction is estimated via an ensemble of simulations over the uncertainty scenarios in $\mathcal S_\theta$:
\begin{align}\label{eq:feasprob}
\widehat{\mathds{P}}[\ve G(\ve d_i,\cdot) \leq \ve 0 \mid \mathcal S_\theta ] := \sum_{(\ve\theta_j,\omega_j)\in\mathcal S_\theta} \mathds{1}[\ve G(\ve d_i,\ve\theta_j)] \: \omega_j
\end{align}
where $\mathds{1}[\cdot]$ stands for the indicator function, such that $\mathds{1}[\ve g] = 1$ if $g_k\added{\leq} 0, \forall k$ and $\mathds{1}[\ve g] = 0$ otherwise; and the samples $\ve\theta_j$ with corresponding weights $\omega_j$ are drawn from the uncertainty distribution $p(\cdot)$. Although effective this approach can prove computationally prohibitive because the number of model simulations grows as the product between the numbers of process parameter samples, $N_d$ and uncertainty scenarios, $N_\theta$. A major challenge for making such methods more efficient thus entails reducing the number of model simulations. The nested sampling algorithm described next aims precisely at sampling the feasibility probability map more uniformly in order to generate a larger number of points with high reliability value. 

\nomenclature[M]{$\omega$}{weight}
\nomenclature[M]{$\mathds{E}[\cdot]$}{expected value}
\nomenclature[M]{$\mathds{1}[\cdot]$}{indicator function}
\nomenclature[M]{$\mathcal{S}$}{sampling set}
\nomenclature[M]{$\mathcal{DS}$}{sampled design space}
\nomenclature[M]{$N$}{number of samples}
\nomenclature[M]{$n$}{dimension}

\begin{algorithm}[tb]
\caption{Standard Monte Carlo Sampling for Design Space Characterization} \label{alg:MC}
\begin{algorithmic}[1]
\State \textbf{input} $\mathcal S_d = \{\ve d_i \in \mathcal K : i=1,\ldots, N_d\}$, $\mathcal S_\theta = \{(\ve \theta_j, \omega_j) \sim p(\cdot) : j=1,\ldots, N_\theta \}$
\State $\mathcal{DS} \gets \emptyset$
\ForAll{$\ve d_i \in \mathcal S_d$}
    \State $\mathcal{DS} \gets \mathcal{DS} \cup \left\{ \left( \ve d_i\: ;\: \widehat{\mathds{P}}[\ve G(\ve d_i,\cdot) \leq \ve 0 \mid \mathcal S_\theta ] \right) \right\}$
\EndFor
\State \Return $\mathcal{DS}$
\end{algorithmic}
\end{algorithm}

\subsection{Nested Sampling for Design Space Characterization}

Nested sampling\cite{skilling2004nested} is a Monte Carlo technique for estimating the evidence in Bayesian parameter estimation. Specifically, given a prior distribution $\pi(\ve \theta)$ and a likelihood function $\mathcal L(\ve\theta)$, the posterior distribution $p(\ve\theta)$ can be inferred from Bayes' theorem:
\begin{align*}
p(\ve\theta) = \frac{\mathcal L(\ve \theta)\:\pi(\ve\theta)}{Z}
\end{align*}
where $Z := \mathds{E}[\mathcal L(\ve\theta) \mid \ve\theta\sim\pi(\cdot)] = \int\!\cdots\!\int_{\mathds{R}^{n_\theta}} \mathcal L(\ve \theta)\:\pi(\ve\theta)\:{\rm d}\ve\theta$ denotes the Bayesian evidence. The algorithm proceeds by progressively sampling in nested contours of increasing likelihood, so as to maintain a dense enough sample in regions of higher likelihood---refer to \replaced{Appendix~A of the Supplementary Information}{the appendix} for further details. This feature is an important motivation for investigating nested sampling in probabilistic DS applications. Before delving into the details of the algorithm, a quick comparison of nested sampling for parameter estimation and DS characterization is presented in Table~\ref{tab:nscompare}.

\nomenclature[M]{$\pi(\cdot)$}{prior distribution}
\nomenclature[M]{$\mathcal{L}(\cdot)$}{likelihood function}
\nomenclature[M]{$Z$}{evidence}

\begin{table}[tb]
\caption{Comparison of nested sampling for parameter estimation and design space characterization.\label{tab:nscompare}}
\begin{tabular}{p{.45\textwidth}|p{.5\textwidth}}
\toprule
\multicolumn{1}{c|}{\bf Bayesian parameter estimation} & \multicolumn{1}{c}{\bf Bayesian design space characterization}\\
\midrule
\rowcolor{lightgray}
\multicolumn{2}{c}{\bf Outcome~~~~~~~}\\
\multicolumn{2}{c}{}\\[-.9em]
An estimate of the Bayesian evidence $Z$ & A set of samples of given size within the probabilistic design space at the desired reliability value $\alpha^*$, $\{\ve d_i\in\mathcal D_{\alpha^*} : i=1,\ldots, N_d\}$ \\[.2em]
A set of weighted samples from the posterior, $\{(\ve\theta_j,\omega_j)\sim p(\cdot):j=1,\ldots, N_\theta\}$ &  A set of samples at a lesser reliability value $\alpha<\alpha^*$, $\{\ve d_k\in\mathcal D_{\alpha_k} : \alpha_k<\alpha^*, k=1,\ldots, N'_d\}$\\
\midrule
\rowcolor{lightgray}
\multicolumn{2}{c}{\bf Challenge~~~~~~~}\\
\multicolumn{2}{c}{}\\[-.9em]
How to efficiently sample the prior distribution subject to a likelihood constraint, $\mathcal L(\ve\theta)>\mathcal L_{\rm min}$ so as to have a high acceptance rate and uniform spread & How to efficiently sample the design space subject to a feasibility probability constraint, $\widehat{\mathds P}[{\bf G}(\ve d,\cdot)\leq 0]>\widehat{\mathds P}_{\rm min}$ so as to have a high acceptance rate and uniform spread\\
\bottomrule
\end{tabular}
\end{table}

\nomenclature[S]{min}{minimal value}
\nomenclature[S]{L}{live points}
\nomenclature[S]{R}{replacement proposals}
\nomenclature[S]{iter}{iterations}
\nomenclature[S]{$*$}{target value}

A pseudo-code of nested sampling for probabilistic DS characterization is presented in Algorithm~\ref{alg:ns_for_ds}. The inputs to this algorithm comprise (line~\ref{alg:ns_ds:line:ini}): the knowledge space, $\mathcal K$; the desired reliability value $\alpha^*$; the number of replacement candidates per iteration $N_{\rm R}$; an initial sample set of design parameters of size $N_{\rm L}$ within the knowledge space $\mathcal S_{\rm L} = \{\ve d_i \in \mathcal K : i=1,\ldots, N_{\rm L}\}$; and set of weighted samples of size $N_\theta$ describing the model parameter uncertainty $\mathcal S_\theta = \{(\ve \theta_j,\omega_j) \sim p(\cdot) : j=1,\ldots, N_\theta \}$. Recall that the latter may either be drawn from a closed-form distribution (with equal weights), e.g. a Gaussian distribution with a given mean and covariance matrix, or be the result of a Bayesian parameter estimation using MCMC or nested sampling again.

\begin{algorithm}[tb]
\caption{Nested Sampling for Design Space Characterization}
\label{alg:ns_for_ds}
\begin{algorithmic}[1]
\State \textbf{input} $\mathcal K$, $\mathcal S_{\rm L} = \{\ve d_i \in \mathcal K : i=1,\ldots, N_{\rm L}\}$, $\mathcal S_\theta = \{(\ve\theta_j,\omega_j) \sim p(\cdot) : j=1,\ldots, N_\theta \}$, $\alpha^*$, $N_{\rm R}$ \label{alg:ns_ds:line:ini}
\State $\mathcal{DS} \gets \emptyset$
\While {$\exists \ve d_i \in \mathcal S_{\rm L} : \widehat{\mathds{P}}[\ve G(\ve d_i,\cdot) \leq \ve 0 \mid \mathcal S_\theta ] < \alpha^*$}
    \State $\widehat{\mathds P}_{\rm min} \gets \min\{ \widehat{\mathds{P}}[\ve G(\ve d_i,\cdot) \leq \ve 0 \mid \mathcal S_\theta ] : d_i \in \mathcal S_{\rm L} \}$\label{alg:ns_ds:line:Pmin}
    \State $\ve d_{\rm min} \gets \arg\min\{ \widehat{\mathds{P}}[\ve G(\ve d_i,\cdot) \leq \ve 0 \mid \mathcal S_\theta ] : d_i \in \mathcal S_{\rm L} \}$
    \State $\mathcal S_{\rm R} \gets \{\ve d_k\in\mathcal K : k=1,\ldots, N_{\rm R} \}$ \label{alg:ns_ds:line:proposal}
    \ForAll {$\ve d_k \in \mathcal S_{\rm R}$}
        \If {$\widehat{\mathds{P}}[\ve G(\ve d_k,\cdot) \leq \ve 0 \mid \mathcal S_\theta ] > \widehat{\mathds{P}}_{\rm min}$} \label{alg:ns_ds:line:test}
            \State $\mathcal S_{\rm L} \gets \mathcal S_{\rm L} \cup \{ \ve d_k \} \setminus \{ \ve d_{\rm min} \}$
            \State $\mathcal{DS} \gets \mathcal{DS} \cup \left\{ \left( \ve d_{\rm min}\: ;\: \widehat{\mathds{P}}_{\rm min} \right) \right\}$
        \EndIf
    \EndFor
\EndWhile
\ForAll{$\ve d_i \in \mathcal S_{\rm L}$}\label{alg:ns_ds:line:loop2.1}
    \State $\mathcal{DS} \gets \mathcal{DS} \cup \left\{ \left( \ve d_i\: ;\: \widehat{\mathds{P}}[\ve G(\ve d_i,\cdot) \leq \ve 0 \mid \mathcal S_\theta ] \right) \right\}$\label{alg:ns_ds:line:loop2.2}
\EndFor
\State \Return $\mathcal{DS}$
\end{algorithmic}
\end{algorithm}

The $N_{\rm L}$ points in the set $\mathcal S_{\rm L}$ are called \textsl{live points} in the nested sampling literature and may be initialized via uniform sampling in the knowledge space $\mathcal{K}$. Nested sampling starts by estimating the feasibility probability at each of these points according to Equation~\eqref{eq:feasprob}, which uses the discretized uncertainty set $\mathcal S_\theta$. Then each iteration of the algorithm proceeds by sampling $N_{\rm R}\geq 1$ new points within an envelope that encloses the current live points and substituting them with the live point having the least feasibility probability in case of improvement. The replaced live points are called \textsl{dead points} and stored with their feasible probability in the result set $\mathcal{DS}$. Notice that the main benefit of using $N_{\rm R}>1$ is in regards of vectorizing the feasibility probability evaluation, although this parallelization capability is not used for the case studies in this paper. 

Key to the efficiency of nested sampling is the ability to generate replacement candidates that have a satisfactory acceptance rate  (line \ref{alg:ns_ds:line:proposal}). An important motivation behind the adaptation of nested sampling for DS characterization is that a variety of techniques have been developed by the Bayesian estimation community over the years to generate such points. This includes: (i) sampling from an enlarged ellipsoid or multiple enlarged ellipsoids enclosing the current live points\cite{mukherjee2006nested,Feroz2009}; or (ii) running a short MCMC from a randomly selected live point\cite{Feroz2009,handley2015polychord}.

Estimating the feasibility probability of the constraints for given process parameter values requires $N_\theta$ process model runs and is by far the most computationally demanding aspect of the algorithm. For $N_{\rm iter}$ nested sampling iterations, the total number of process model runs amounts to $N_\theta (N_{\rm L} + N_{\rm iter} N_{\rm R})$. An acceleration strategy consists in interrupting the evaluation of the feasibility probability $\widehat{\mathds{P}}[\ve G(\ve d_k,\cdot) \leq \ve 0 \mid \mathcal S_\theta ]$ of a replacement candidate $\ve d_k$ (line \ref{alg:ns_ds:line:test}) if the feasibility probability mass accumulated so far and the mass of the remaining samples add up to less than $\hat{\mathds{P}}_{\rm min}$ (line \ref{alg:ns_ds:line:Pmin})---e.g. by summing in decreasing order of probability weights $\omega_j$ in Equation~\eqref{eq:feasprob} for maximal efficiency. Such an interruption is indeed a guarantee that the replacement candidate should be rejected. Of course, the consequence is that the feasibility probability $\widehat{\mathds{P}}[\ve G(\ve d_k,\cdot) \leq \ve 0 \mid \mathcal S_\theta ]$ of any rejected point will not be estimated accurately when this acceleration strategy is used. This is the main reason why we do not append the rejected points to the result set $\mathcal DS$ in Algorithm~\ref{alg:ns_for_ds}, which is in agreement with our primary objective of uncovering $N_{\rm L}$ points with a feasibility probability greater than $\alpha^*$.

The main iteration terminates when the feasibility probability of all the live points is no less than the target reliability value $\alpha^*$. All of these points are appended with their corresponding feasibility probabilities to the result set $\mathcal{DS}$ (line \ref{alg:ns_ds:line:loop2.2}) that already contains the dead points. It may happen that the number of points with a feasibility probability greater than $\alpha^*$ in $\mathcal{DS}$ is larger than $N_{\rm L}$ when multiple live point replacements occur during the final iteration. Also, Algorithm~\ref{alg:ns_for_ds} will terminate prematurely in the event that $\mathcal{D}_{\alpha^*}$ is empty, which occurs when the model uncertainty is too large or the reliability value $\alpha^*$ is too high. This situation can nonetheless be detected, e.g. by monitoring the rate of improvement in the feasibility probability over a given number of iterations, and the algorithm can be stopped when progress is too slow.

An implementation of Algorithm~\ref{alg:ns_for_ds} in a Python package called {\sf DEUS} (standing for \textsl{DEsign under Uncertainty using Sampling techniques}) is available as part of the Supplementary Material to the paper. At each iteration the replacement candidates are generated by sampling in a single ellipsoid enclosing the current live points. Accordingly, the tuning parameters for a probabilistic DS computation in {\sf DEUS} include: (i) the number of live points, $N_{\rm L}$; (ii) the number of replacement candidates at each iteration, $N_{\rm R}$; (iii) the initial enlargement factor of the ellipsoid; and (iv) the shrinking rate of that enlargement factor at each iteration. This package is used to solve all of the numerical case studies presented below, using as default parameters an initial enlargement factor of 30\% and a shrinking rate of 0.2 for the ellipsoids\cite{Feroz2009} as well as $N_{\rm R}=8$ replacement proposals. Despite the fact that no parallelization is implemented in {\sf DEUS} at present, choosing an $N_{\rm R}$ value greater than $1$ can reduce the overhead caused by reconstructing an ellipsoid around the live points too frequently. Finally, an appropriate $N_{\rm L}$ value is highly dependent on the volume ratio between target design space and knowledge space, with a larger volume ratio calling for a larger $N_{\rm L}$ to maintain a suitable sample density across the DS. In our experience it is advisable to start with a low $N_{\rm L}$ to establish that $\mathcal{D}_{\alpha^*}$ is not empty, then refine the characterization by increasing $N_{\rm L}$. \added{Like any other sampling-based approach, there is no guarantee that {\sf DEUS} will not miss part of the feasible region. But the likelihood of this happening can be reduced by either increasing the enlargement factor of the ellipsoids, decreasing the related shrinking rate, or increasing the number of replacement proposals per iteration.}

\nomenclature[A]{{\sf DEUS}}{DEsign under Uncertainty using Sampling techniques}

\subsection{Exploitation of the Results}

Sampling-based approaches to probabilistic DS characterization return a set of points and their corresponding feasibility probabilities. Turning these results into a format that can be exploited by the practitioner is clearly important. It is straightforward to represent a probabilistic design space in one or two dimensions graphically, for instance using a range of colors to depict different feasibility probability levels. This representation may also be used in three or four dimensions by splitting the range of the extra dimensions onto a trellis chart---see for instance the Suzuki reaction case study below.

But, more generally, a graphical presentation might be insufficient for the practitioner to easily test the feasibility probability of a given set of design parameters or to convey the resulting design space to third-parties in a concise way. Alternative representations include:
\begin{itemize}
\item Inscribing a simple shape---e.g. a box or an ellipsoid---within the envelope defined by the samples above a given reliability level in the data set $\mathcal{DS}$ returned by Algorithm~\ref{alg:MC} or \ref{alg:ns_for_ds}. Under the assumption that the samples describe a convex set, one may proceed by first constructing a polyhedral representation of the set then fitting the desired shape within that polyhedron. The latter problem may be solved efficiently using convex optimization, but constructing a set of hyperplanes to describe the convex hull of a set of points can be computationally burdensome (NP-hard)\cite{Boyd2004}. This approach is furthermore inadequate for feasibility probability maps presenting nonconvex iso-contours or multiple modes.
\item Fitting a non-parametric model to approximate the full feasibility probability map across the knowledge space. For instance, a a multilayer perceptron\cite{hastie2009,Goodfellow2016} (MLP)---a class of feed-forward artificial neural network with multiple (hidden) layers and non-linear activation functions---can be trained on the labelled data set $\mathcal{DS}$. The trained MLP provides a computationally cheap surrogate, transforming the DS samples into a suitable form for exploitation. It may be used readily to predict the feasibility probability for any design parameters $\ve d\in\mathcal K$. It may also be embedded into a design-centering optimization problem\cite{harwood2017} for finding a subset of points, e.g. in the form of a box or an ellipsoid, with feasibility probability above a desired reliability value. This is supported by recent advances in deterministic global optimization with neural networks embedded\cite{Schweidtmann2019}.
\item Fitting a multinomial classifier\cite{hastie2009} to separate the knowledge space into two or more subregions corresponding to different feasibility probability ranges---e.g. below and above the reliability value $\alpha^*$. Similar to non-parametric regression an MLP has the ability to distinguish data sets that are not linearly separable and may be trained on the labelled data set $\mathcal{DS}$. The trained classifier then provides a cheap way of estimating the probability range of any design parameters $\ve d\in\mathcal K$ and a softmax function may be used in the final layer of the MLP to estimate the probability of $\ve d\in\mathcal K$ belonging to a given range. Like a non-parametric regression model such a classifier may also be embedded into a design-centering problem.
\end{itemize}

\nomenclature[A]{MLP}{Multi-Layer Perceptron}

Our focus in the remainder of the paper is on MLP to approximate the full feasibility probability map. The MLP of interest has a single neuron in its output layer---whose state represents the feasibility probability---and its input layer comprises exactly $n_d$ neurons---one for each process parameter. A key, yet arduous, decision is selecting the numbers of hidden layers and hidden neurons in the MLP. At least one hidden layer is necessary since the feasibility probability map is generally nonlinear. An MLP with a single hidden layer is capable of approximating any continuous functions under mild assumptions on the activation function but it can take an arbitrary large number of neurons in that hidden layer to meet a desired accuracy\cite{Hornik1991}. Deep learning mitigates this problem by including additional hidden layers. However, it is important to keep in mind that the feasibility probability is inherently noisy due to the model uncertainty discretization. Therefore, the MLP may also be called upon to play a role in filtering this noise and avoiding overfitting of the data by restricting its number of hidden layers or neurons. Since the numerical case studies below have a handful of process parameters we consider hidden layers with a few dozen neurons as a rule of thumb. A more systematic analysis is beyond the scope of this paper.

\subsection{Illustrative Example (Continued)}

\begin{figure}[tbp]
\centering
\includegraphics[width=.98\textwidth]{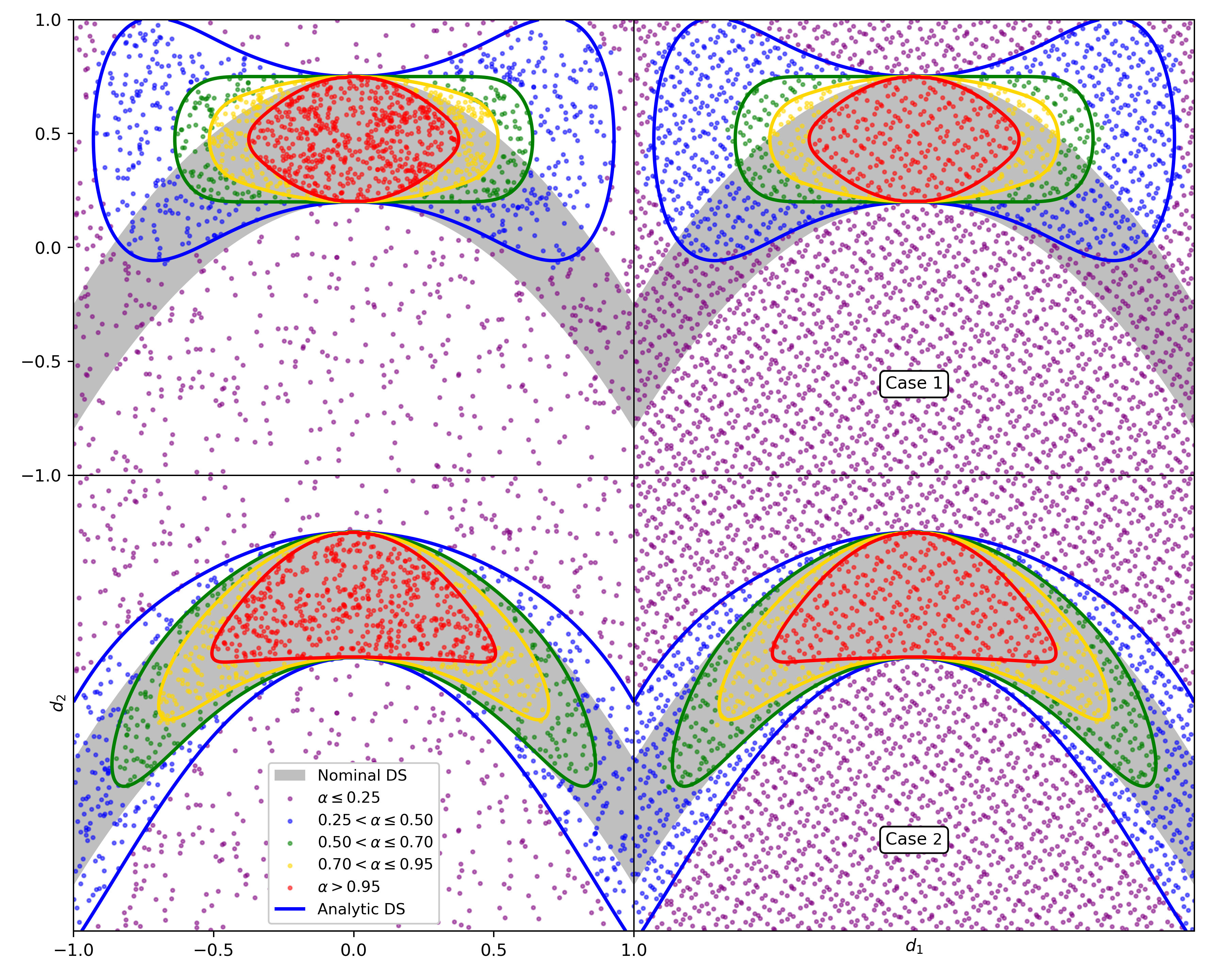}
\caption{Comparison of sampling-based techniques to determine the probabilistic design space in Equation~\eqref{eq:prob_ds_example}. The top and bottom plots are for the uncertainty scenarios $\theta \sim \mathcal N(0,1)$ and $\theta \sim \mathcal N(1,\sqrt{0.3})$, respectively. The left plots show the results of nested sampling using Algorithm~\ref{alg:ns_for_ds}. The right plots show the results of the standard Monte Carlo in Algorithm~\ref{alg:MC} using quasi-random Sobol sampling. The nominal design space (Equation \ref{eq:nom_ds_example}) is represented with the gray-shaded area. Iso-reliability contours for the actual probabilistic design space at different reliability values $\alpha$ (Equation~\ref{eq:prob_ds_example}) are shown in solid lines with different colors.}
\label{fig:sample_DS_example}
\end{figure}

We revisit our simple case study concerned with the determination of the design space defined by the CQA constraint \eqref{eq:CQA_ds_example} using the uncertain model \eqref{eq:model_ds_example}. The left plots in Figure~\ref{fig:sample_DS_example} are the results of nested sampling (Algorithm~\ref{alg:ns_for_ds}) with a reliability target of $\alpha^*=0.95$, $N_{\rm L}=500$ live points, and $N_\theta=100$ samples drawn from either $\mathcal N(0,1)$ (top plot) or $\mathcal N(1,\sqrt{0.3})$ (bottom plot). Overall, the resulting samples---depicted in a different color according to reliability membership---are in excellent agreement with the exact subregions delimited by the solid contours (Equation~\ref{eq:prob_ds_example}). Nevertheless a few of the points are misclassified due to the approximation of the feasibility probability based on Equation~\eqref{eq:feasprob} which uses a finite number of uncertainty samples. This discrepancy can be removed by increasing the number of uncertainty samples $N_\theta$, yet at the cost of a higher computational burden since the total number of model evaluations is directly proportional to $N_\theta$.

The right plots in Figure~\ref{fig:sample_DS_example} are generated using the Monte Carlo Algorithm~\ref{alg:MC} with Sobol sampling. These results are again in excellent agreement with the exact subregions, despite a few samples being misclassified for the same reason as earlier with nested sampling since the same approach is used to estimate the feasibility probability. For a quantitative comparison between the two algorithms, notice that the total number of model evaluations is identical in both cases (322,000 and 369,200 model evaluations under the uncertainty scenarios $\theta \sim \mathcal N(0,1)$ and $\theta \sim \mathcal N(1,\sqrt{0.3})$, respectively) and no acceleration strategy is applied for nested sampling. A breakdown of the samples generated with Algorithms~\ref{alg:MC} and~\ref{alg:ns_for_ds} within different feasibility probability ranges is furthermore reported in Table~\ref{tab:sample_DS_example}. Recall that unsuccessful proposals in the nested sampling approach are currently discarded, which is the reason why the total number of samples is lower with this approach compared to Monte Carlo sampling. Regardless, nested sampling enables a much denser sampling of the targeted reliability region, as controlled by the number of live points and the gradual sampling towards higher reliability regions. It would take three to four times more samples---and hence model evaluations---for the standard Monte Carlo approach to achieve a similar concentration of points within the desired reliability region of $0.95$ in this case study. 

\begin{table}[tb]
\caption{Comparison between the samples generated via standard Monte Carlo with Sobol sampling (Algorithm~\ref{alg:MC}) and nested sampling (Algorithm~\ref{alg:ns_for_ds}) within different reliability ranges.} 
\label{tab:sample_DS_example}
\begin{tabular}{lrr}
\toprule
{\bf Reliability} & {\bf Samples drawn by standard} & {\bf Samples drawn by nested}\\
{\bf value} & \multicolumn{1}{l}{\bf Monte Carlo (Algorithm~\ref{alg:MC})} & \multicolumn{1}{l}{\bf sampling (Algorithm~\ref{alg:ns_for_ds})}\\
\midrule
\rowcolor{lightgray}
& \multicolumn{2}{c}{\bf Uncertainty scenario: $\theta \sim \mathcal N(0,1)$}\\
& \multicolumn{2}{c}{}\\[-.9em]
$0.95 \leq \alpha \hphantom{< 1.00}$    & \replaced{247}{267}   & 500 \\
$0.70 \leq \alpha < 0.95$               & \replaced{139}{156}   & \replaced{231}{281} \\
$0.50 \leq \alpha < 0.70$               & \replaced{172}{212}   & \replaced{181}{172} \\
$0.25 \leq \alpha < 0.50$               & \replaced{666}{654}   & \replaced{418}{388} \\
$\hphantom{0.00 \leq}~\alpha < 0.25$    & \replaced{2,025}{1,960} & \replaced{532}{524} \\
\midrule
Total & 3,249 & \replaced{1,862}{1,865} \\
\midrule
\rowcolor{lightgray}
& \multicolumn{2}{c}{\bf Uncertainty scenario: $\theta \sim \mathcal N(1,\sqrt{0.3})$}\\
& \multicolumn{2}{c}{}\\[-.9em]
$0.95 \leq \alpha \hphantom{< 1.00}$    & \replaced{331}{354}   & \replaced{500}{509} \\
$0.70 \leq \alpha < 0.95$               & \replaced{194}{225}   & \replaced{210}{259} \\
$0.50 \leq \alpha < 0.70$               & \replaced{212}{256}   & \replaced{200}{171} \\
$0.25 \leq \alpha < 0.50$               & \replaced{490}{345}   & \replaced{321}{204} \\
$\hphantom{0.00 \leq}~\alpha < 0.25$    & \replaced{2,022}{2,069}  & \replaced{572}{597} \\
\midrule
Total & 3,249 & \replaced{1,803}{1,740} \\
\bottomrule
\end{tabular}
\end{table}

\begin{figure}[tbp]
\centering
\includegraphics[width=.45\textwidth]{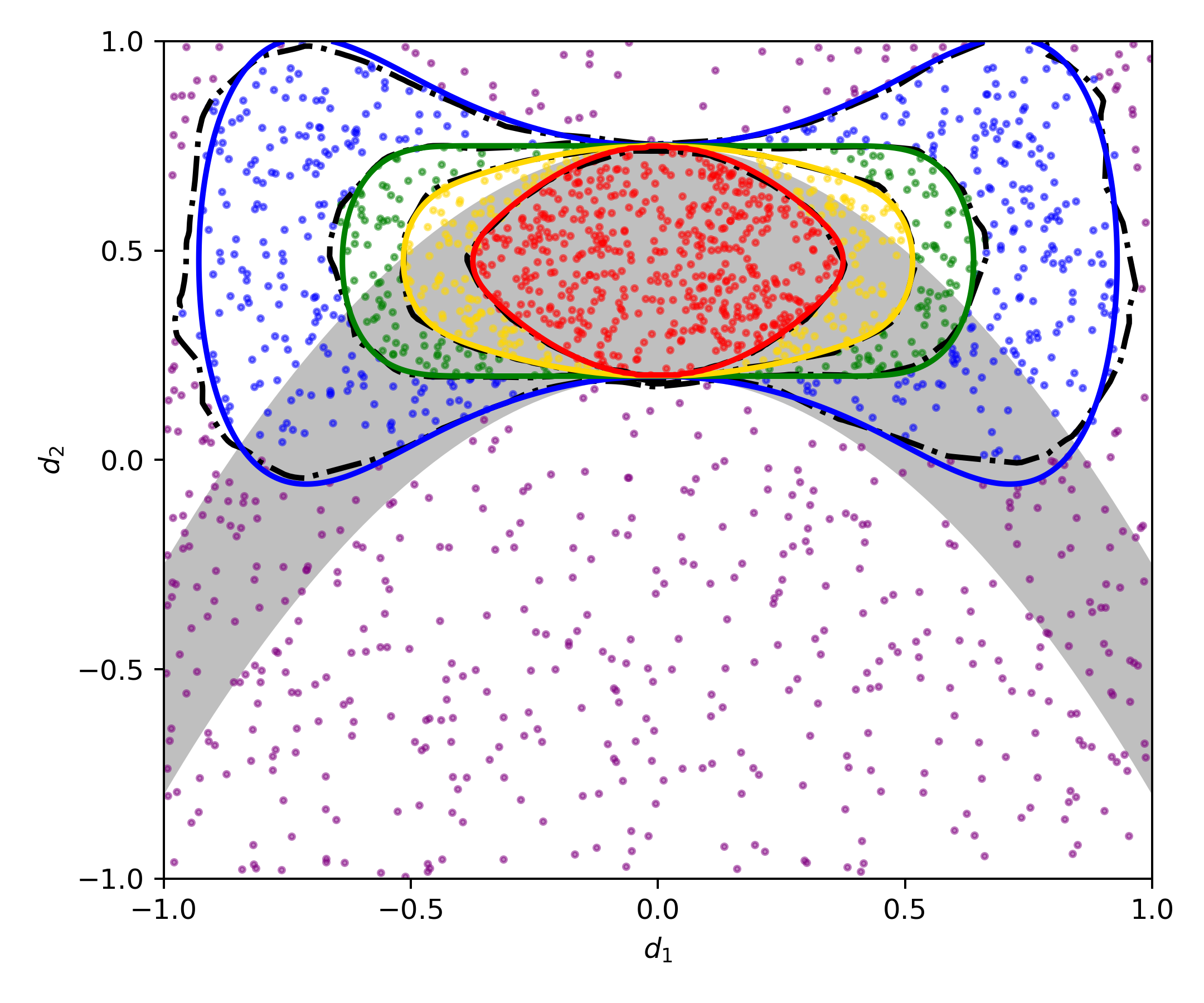}\hfill
\includegraphics[width=.51\textwidth]{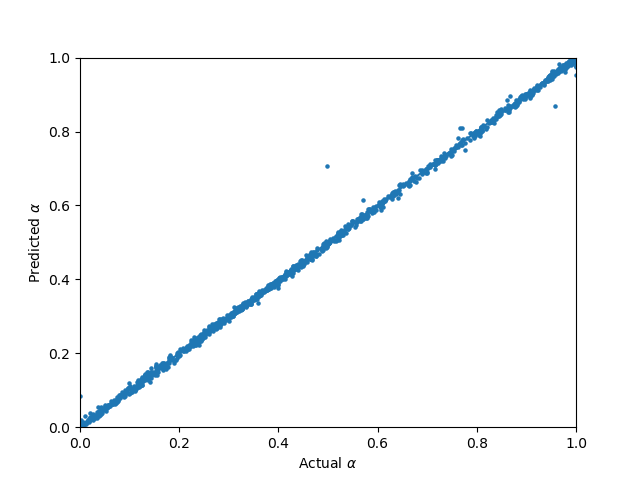}
\caption{Exploitation of the nested sampling results using an MLP to fit the full feasibility probability map in the uncertainty scenario $\theta \sim \mathcal N(0,1)$. The left plot compares the contours of the fitted MLP (dashed lines) with the actual DS (solid color lines) at different reliability values $\alpha$. The gray-shaded area depicts the nominal design space (Equation~\ref{eq:nom_ds_example}). The right plot shows the parity plot of the trained MLP.}
\label{fig:MLPfit_DS_example}
\end{figure}

The left plot in Figure~\ref{fig:MLPfit_DS_example} compares several iso-contours of the exact feasibility probability map (Equation~\ref{eq:nom_ds_example}) with those of an MLP trained on the nested sampling results under the uncertainty scenario $\theta \sim \mathcal N(0,1)$. A good agreement is obtained in the upper reliability region where the concentration of samples is high, while the MLP predictions somewhat deviate in the lower reliability region where the sample points are scarcer---compare Table~\ref{tab:sample_DS_example}. The MLP under consideration here has four hidden layers, with 16, 32, 32 and 16 artificial neurons, respectively. Its training was conducted in the Python package {\sf Keras} \cite{chollet2015keras} that is implemented on top of the package {{\sf TensorFlow}\cite{tensorflow2015-whitepaper}}, using mini-batch gradient descent. The resulting parity plot is shown on the right plot in Figure~\ref{fig:MLPfit_DS_example}. The regression quality is generally good but for a few of the points that present a rather large deviation from the training set. Part of this mismatch is attributed to the fact that the data points in the training set itself are noisy due to the uncertainty discretization. Moreover, the actual feasibility probability map presents discontinuities at the points $\ve d= (0, 0.20)^\intercal$ and $\ve d= (0, 0.75)^\intercal$ so the MLP has difficulties capturing the very stiff probability variations near such points.

\section{Industrial Case Studies}

We demonstrate the applicability of the proposed nested sampling approach on two case studies of industrial relevance. The first one considers a Michael addition reaction in a continuous reactor and enables a comparison with recent work on flexibility-based algorithms\cite{laky2019optimization} for probabilistic DS characterization. The second one investigates a biphasic Suzuki coupling reaction performed in fully batch mode\cite{garcia2015definition} and gives the opportunity to determine a four-dimensional DS with a more complex mechanistic process model. 


\subsection{Michael Addition Reaction}

This case study considers the following Michael addition reaction in a continuous stirred-tank reactor (CSTR):\cite{laky2019optimization}
\begin{align*}
{\sf AH} + {\sf B} & \ce{->[$k_1$]} {\sf A^-} + {\sf BH^+} \\
{\sf A^-} + {\sf C} & \ce{<=>[$k_2$][$k_3$]} {\sf AC^-} \\
{\sf AC^-} + {\sf AH} & \ce{->[$k_4$]} {\sf A^-} + {\sf P} \\
{\sf AC^-} + {\sf BH^+} & \ce{->[$k_5$]} {\sf P} + {\sf B}
\end{align*}
where {\sf AH} is the Michael donor; {\sf C} is the Michael acceptor; {\sf B} is a base; {\sf BH$^+$}, {\sf A$^-$}, and {\sf AC$^-$} are reaction intermediates; and {\sf P} is the product. A complete statement of the steady-state process model can be found in Appendix~B \added{of the Supplementary Information}. All of the kinetic constants are considered uncertain in this model, following a multivariate normal distribution $(k_1,\ldots,k_5)\sim\mathcal N(\mu_{\ve k},\Sigma_{\ve k})$.

\nomenclature[A]{CSTR}{Continuous Stirred-Tank Reactor}
\nomenclature[M]{$k$}{kinetic rate constant}
\nomenclature[M]{$\uptau$}{residence time [min]}
\nomenclature[M]{$R$}{molar ratio [mol\,mol$^{-1}$]}
\nomenclature[M]{$T$}{temperature [\textdegree{C}]}
\nomenclature[M]{$y$}{gas molar fraction [ppm]}

The knowledge space is defined in terms of two process parameters: (i) the molar ratio between the concentration of {\sf AH} and {\sf B} in the feed, $R_{\sf AH\mid B}\in[10,30]$; and (ii) the residence time in the CSTR, $\uptau\in[400,1400]~\rm (min)$. The feasible operating region is furthermore limited by two CQA constraints: (i) conversion of feed {\sf C} greater than 90\%; and (ii) residual concentration of {\sf AC$^-$} smaller than $2\,\rm mmol\,L^{-1}$.  

\begin{figure}[tbp]
\centering
\includegraphics[width=.9\linewidth]{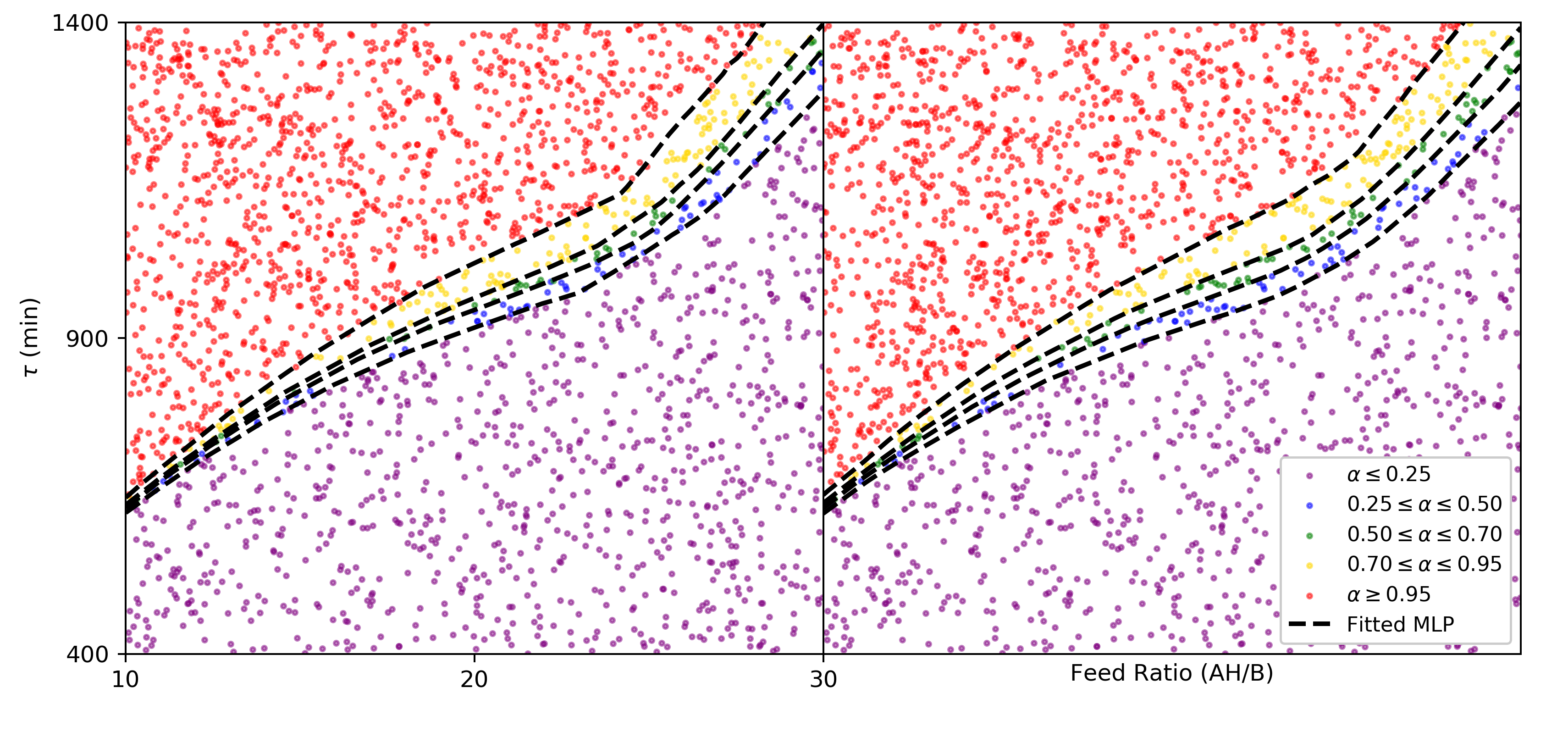}
\caption{Comparison of probabilistic DS for the Michael addition reaction computed using nested sampling (Algorithm~\ref{alg:ns_for_ds}) with $N_\theta=100$ (left) and $1,000$ (right) uncertainty scenarios. Samples from the probabilistic design space belonging to different reliability ranges are shown in different colors. Iso-reliability contours from the fitted MLP are shown in black dashed lines.}
\label{fig:comp_DS_michael}
\end{figure}

Like in the illustrative example above, we apply nested sampling (Algorithm~\ref{alg:ns_for_ds}) with a reliability target of $\alpha^*=0.95$, now considering $N_{\rm L}=1,000$ live points to enable a high density of points. The resulting sample set ($\mathcal{DS}$) is then regressed using an MLP comprising 4~hidden layers (with 16, 32, 32 and 16 neurons in each layer) to approximate the full feasibility probability map. The two plots in Figure~\ref{fig:comp_DS_michael} are for $N_\theta=100$ and $1,000$ uncertainty samples. The corresponding computational statistics are reported in Table~\ref{tab:comp_DS_michael}, both without and with the acceleration strategy. The effect of increasing the number of uncertainty scenarios by 10 fold is rather small in Figure~\ref{fig:comp_DS_michael}, while the increase in CPU time is naturally close to a factor a 10 due to the corresponding increase in the number of function evaluations. Though differences can be noted between both fitted MLPs, with the MLP on the right plot ($N_\theta=1,000$) predicting smoother iso-reliability contours due to a lower noise in the estimated feasibility probabilities compared to the MLP on the left plot ($N_\theta=100$). It is also worth noting that the effect of the acceleration---about 5\% reduction in the number of function evaluations and $10\%$ reduction in CPU time---is rather modest. This is attributed to a relatively low fraction of rejected proposals due to the shape of the DS. 

\nomenclature[A]{CPU}{Central Processing Unit}

\begin{table}[tb]
\caption{Computational statistics for the Michael addition reaction using nested sampling (Algorithm~\ref{alg:ns_for_ds}) with different model uncertainty representations.} 
\label{tab:comp_DS_michael}
\begin{threeparttable}
\begin{tabular}{rrrrr}
\toprule
{\bf Uncertainty} & \multicolumn{2}{c}{\bf Nested sampling} & \multicolumn{2}{c}{\bf Nested sampling}\\
{\bf scenarios,} & \multicolumn{2}{c}{\bf with acceleration} & \multicolumn{2}{c}{\bf without acceleration}\\
\cmidrule{2-5}
{\bf $N_\theta$} & \# model eval. & CPU sec.\tnote{1} & \# model eval. & CPU sec.\tnote{1} \\
\midrule
100   & 226,884 & 114 & 236,000 & 128 \\
1,000 & 2,287,174 & 1,070 & 2,400,000 & 1,207 \\
\bottomrule
\end{tabular}
\begin{tablenotes}
\item[1] CPU times obtained on a single core of \added{{\sf AMD Ryzen 5 2600X}} processor.
\end{tablenotes}
\end{threeparttable}
\end{table}

Finally, we compare the nested sampling results to the results obtained by \citet{laky2019optimization} using an optimization approach based on flexibility analysis. It is worth reiterating that the algorithms therefrom determine or approximate the set-membership counterpart $\widehat{\mathcal D}_\alpha$ (Equation~\ref{eq:robust_ds_definition}) to the probabilistic design space ${\mathcal D}_\alpha$, which happens to be quite conservative\cite{laky2019optimization}. These optimization-based methods are consistently faster than nested sampling (by a factor of 10 or more) when solving the optimization problems using a (local) gradient-based algorithm. But nested sampling is nevertheless competitive against these methods when a global optimization solver is used.

\subsection{Suzuki Coupling Reaction}

This case study investigates the Suzuki coupling reaction between a boronic ester ({\sf SM1}) and an organohalide ({\sf SM2}) to produce a desired pharmaceutical intermediate ({\sf P1}) and a dimeric impurity ({\sf Imp1}) related to {\sf SM1} \cite{garcia2015definition}. The reaction is biphasic and conducted in batch mode. The gaseous phase consists of an inert gas with traces of $\rm O_2$; the liquid phase consists of a mixture of water and tetrahydrofuran (THF) as solvent for 17 chemical species that participate in 12 reactions---3 of which are reversible and 1 is considered instantaneous. A mechanistic, kinetic-based model is available describing the changes in composition of the liquid and gas phases during the batch---see Appendix~C \added{of the Supplementary Information} for a complete statement. The temperature-dependent reaction rates are modeled with Arrhenius equations and kinetic parameters that were verified experimentally\cite{garcia2015definition}. The pre-exponential factors of all 14 kinetically-limited reactions are considered to be uncertain here, following normal distributions with standard-deviations equal to 15\% of the nominal values. 

The reactor has a large number of design and operation parameters so we consider a reduced DS characterization problem in terms of four parameters only: (i) the batch duration, $\uptau\in[75,300]~\rm (min)$; (ii) the equivalent of catalyst, $R_{\sf Pd\mid SM2}\in[0.001,0.003]~\rm (mol~mol^{-1})$; (iii)~the reactor temperature, $T\in[22,64]~\rm (\text{\textdegree{C}})$; and (iv) the molar fraction of \ce{O2} in reactor's head space $y_{\ce{O2}}\in[10,250]~\rm (ppm)$. Moreover, two CQAs limit the feasible operating region at the end of the batch: (i) maximum amount of unreacted {\sf SM2} of $0.001~\rm mol\:mol^{-1}$ for the reaction to be considered complete; and (ii) maximum level of impurity {\sf Imp1} below $0.0015~\rm mol\:mol^{-1}$ for the batch product to be downstream processable.

\begin{figure}[tbp]
\centering
\includegraphics[width=\linewidth]{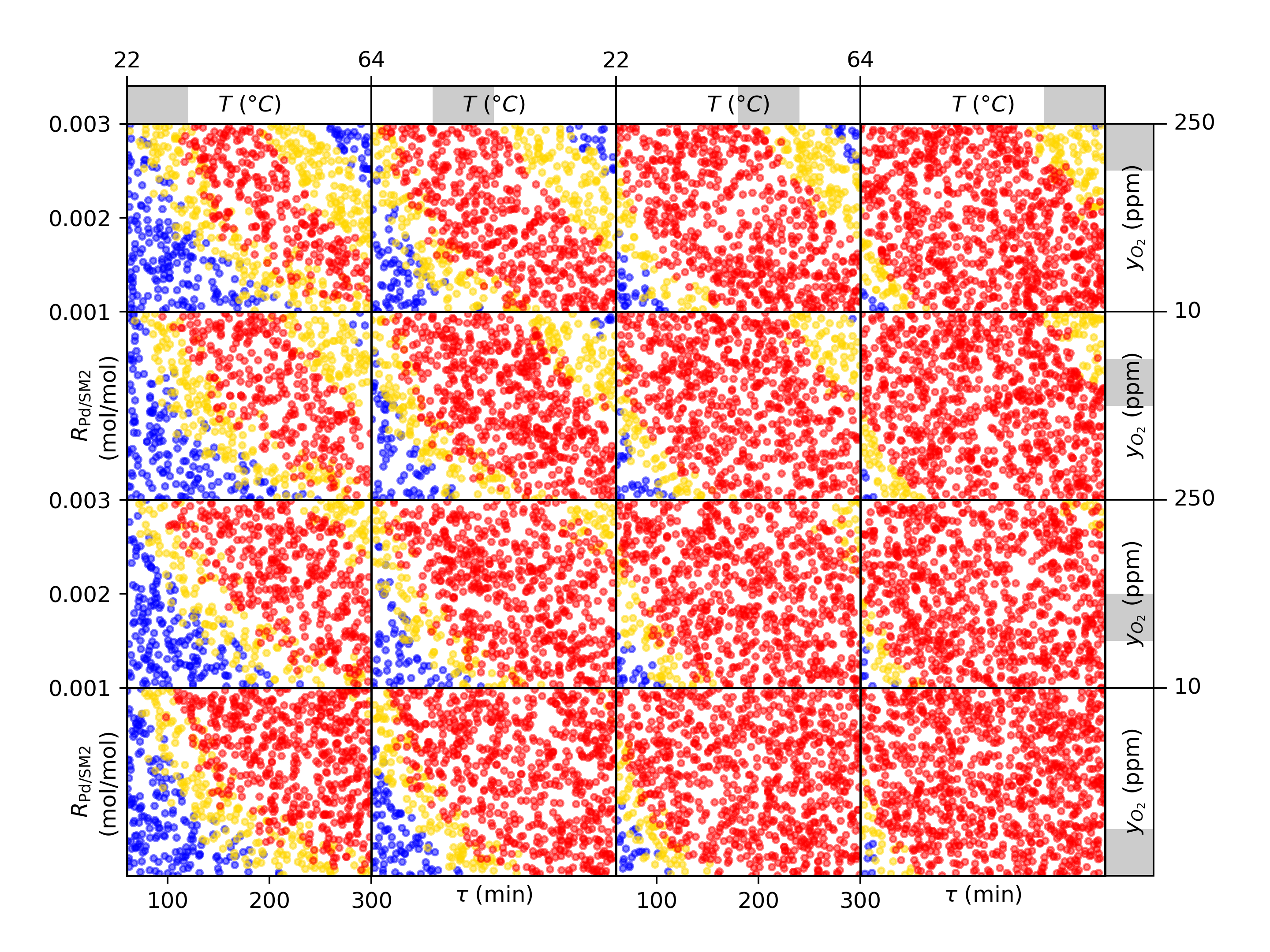}
\caption{Probabilistic DS for the Suzuki coupling reaction computed using nested sampling (Algorithm~\ref{alg:ns_for_ds}) with \replaced{$N_L=10,000$ live points and $N_\theta=1,000$}{$N_\theta=200$} uncertainty scenarios. The outer axes on the trellis chart correspond to the oxygen fraction in head space ($y_{\sf O_2}$) and reactor temperature ($T$); the inner axes on each subplot, to the batch duration ($\uptau$) and catalyst equivalent ($R_{\sf Pd\mid SM2}$). Samples from the probabilistic design space belonging to different reliability ranges are shown in different colors (red: $\alpha\geq 0.85$; yellow: $0.05\leq\alpha<0.85$; blue: $\alpha<0.05$).}
\label{fig:prob_DS_suzuki}
\end{figure}

In order to determine the probabilistic DS we apply nested sampling (Algorithm~\ref{alg:ns_for_ds}) with a reliability target of 
\replaced{$\alpha^*=0.85$}{$\alpha^*=0.95$}. Because the DS now comprises four process parameters \deleted{we increase }the numbers of live points \replaced{needs increasing to ensure a high density of points, so we use}{to} $N_{\rm L}=5,000$\added{ and $10,000$ next}. \added{We furthermore consider two uncertainty descriptions with $N_\theta=200$ and $1,000$ samples, respectively.} For the DS representation we use trellis charts where the ranges of oxygen concentration and temperature are split into four-by-four intervals---indicated by grey bars on the two outer axes; then each subplot is a projection of the points that belong to the particular intervals of oxygen fraction in head space $y_{\sf O_2}$ and temperature $T$ on the plane defined by the batch duration $\uptau$ and catalyst equivalent $R_{\sf Pd\mid SM2}$---the two inner axes. For instance, the chart on Figure~\ref{fig:prob_DS_suzuki} is a representation of the probabilistic DS computed with \added{$10,000$ live points and} \replaced{$1,000$}{$N_\theta=200$} uncertainty scenarios. Notice in particular the effect of the process parameters on the chart, whereby the DS expands upon increasing temperature at constant oxygen level and upon increasing oxygen fraction in head space at constant temperature. Computational statistics for the probabilistic DS are shown in Table~\ref{tab:comp_DS_suzuki}\replaced{.}{ using nested sampling with acceleration---} Overall, it takes \replaced{over 4 days}{over one day of CPU time} to determine the probabilistic DS \added{with $10,000$ live points and $1,000$ uncertainty scenarios. A 12-fold decrease in the CPU time (to about 9 hours) is observed upon reducing to $5,000$ live points and $200$ uncertainty scenarios.} The number of function evaluations is only reduced by about 3\% in applying the acceleration strategy since the fraction of rejected proposals is again low. \added{The corresponding reduction in CPU time is between 4--10\%.} A vectorized implementation of Algorithm~\ref{alg:ns_for_ds} that would compute simulation ensembles or replacement proposals in parallel, could readily reduce the overall time needed and will be the focus of future work. 

\begin{table}[tb]
\caption{Computational statistics for the Suzuki coupling reaction using nested sampling (Algorithm~\ref{alg:ns_for_ds}) with different model uncertainty representations \added{and live points}.} 
\label{tab:comp_DS_suzuki}
\begin{threeparttable}
\begin{tabular}{rrrrrr}
\toprule
{\bf Live} & {\bf Uncertainty} & \multicolumn{2}{c}{\bf Nested sampling} & \multicolumn{2}{c}{\bf Nested sampling}\\
{\bf Points,} & {\bf scenarios,} & \multicolumn{2}{c}{\bf with acceleration} & \multicolumn{2}{c}{\bf without acceleration}\\
\cmidrule{3-6}
{\bf $N_L$} & {\bf $N_\theta$} & \# model eval. & CPU min.\tnote{1} & \# model eval. & CPU min.\tnote{1} \\
\midrule
5,000   & 200   & \added{1,412,390}  & \added{525} & \added{1,451,200} & \added{574}\hphantom{.5} \\
\added{10,000} & \added{1,000} & \added{14,132,975} & \added{6,464} & \added{14,544,000} & \added{6,741}\hphantom{.5} \\
5,000 & 1 & -- & -- & 6,392 & 3.5 \\
\added{10,000} & \added{1} & -- & -- & \added{12,680} & \added{6.8} \\
\bottomrule
\end{tabular}
\begin{tablenotes}
\item[1] CPU times obtained on a single core of \added{{\sf AMD Ryzen 5 2600X}} processor.
\end{tablenotes}
\end{threeparttable}
\end{table}

\begin{figure}[tbp]
\centering
\includegraphics[width=\linewidth]{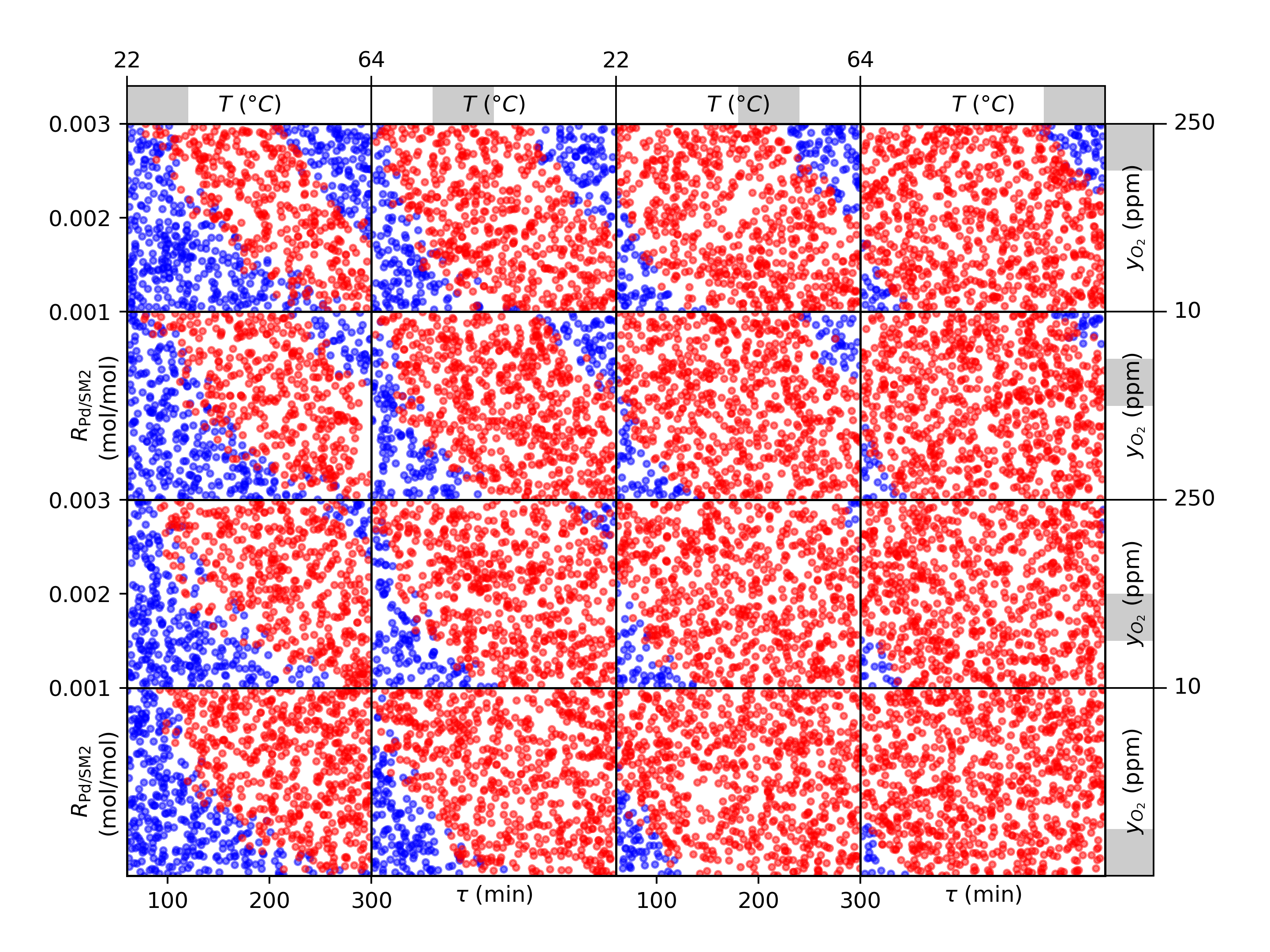}
\caption{Nominal DS for the Suzuki coupling reaction computed using nested sampling (Algorithm~\ref{alg:ns_for_ds})\added{ with $N_L=10,000$ live points}. The outer axes on the trellis chart correspond to the oxygen fraction in head space ($y_{\sf O_2}$) and reactor temperature ($T$); the inner axes on each subplot, to the batch duration ($\uptau$) and catalyst equivalent ($R_{\sf Pd\mid SM2}$). Samples within the nominal design space are depicted in red and those outside are in blue.}
\label{fig:nom_DS_suzuki}
\end{figure}

For comparison, we also compute the nominal DS (Equation~\ref{eq:nom_ds_definition}) by applying nested sampling (Algorithm~\ref{alg:ns_for_ds}), in this case with a single uncertainty scenario ($N_\theta=1$) corresponding to the nominal kinetic parameters and setting the reliability target to $\alpha^*=1$. The chart in Figure~\ref{fig:nom_DS_suzuki} is a representation of this nominal DS\added{ with 10,000 live points}. The \replaced{results}{reported computational times} in Table~\ref{tab:comp_DS_suzuki} \replaced{for 5,000 and 10,000 live points show a dramatic reduction in CPU time}{is  dramatically reduced} to just a few minutes, which is clearly attributed to the fact that a single simulation is needed to assess whether or not a point belongs to the nominal DS instead of a large ensemble of simulations. However, the results in Figures~\ref{fig:prob_DS_suzuki} and~\ref{fig:nom_DS_suzuki} differ significantly as certain points within the nominal DS have a feasibility probability as low as 10--20\%---an extremely poor reliability level in practice. Despite their popularity among practitioners approaches based on nominal parameter values or mean responses fail to quantify reliability and risk and, therefore, their results can be misleading \cite{garcia2015definition,peterson2017}.

\section{Concluding Remarks}

Design space is a key concept in pharmaceutical QbD, helping practitioners develop a better understanding of their manufacturing processes and enhancing regulatory flexibility. In this context, it is of paramount importance to develop computational methods and tools that can provide a clear quantitative representation of the DS---in agreement with the ICH Q8 guideline. Our main focus throughout the paper has been on Bayesian approaches to DS characterization, which determine a feasibility probability that can be used as a measure of reliability and risk.

We have contributed a nested sampling algorithm tailored to the characterization of a probabilistic DS. A hallmark feature of nested sampling is its ability to maintain a given set of live points through regions with increasing probability feasibility until reaching a desired reliability level (or showing it cannot be reached). The algorithm furthermore leverages efficient strategies from Bayesian parameter estimation for generating replacement proposals during the search. Through a simple illustrative example and the case study of a Michael addition reaction, we have established that nested sampling can outperform standard Monte Carlo sampling and be competitive with optimization methods relying on process flexibility concepts, even in low-dimensional DS characterization problems. We have also showcased the use of machine learning techniques to reconstruct a feasibility probability map based on the sampled design space, which can be more easily exploited by the practitioner. In the second case study of a Suzuki coupling reaction we have shown that nested sampling is effective for larger DS characterization with a handful of process parameters, in the presence of a complex dynamic model and realistic model uncertainty---a class of problems currently out of reach for flexibility-based optimization techniques. 

A major impediment facing nested sampling---and other sampling-based techniques---for probabilistic DS characterization in higher-dimensional problems, for instance in multi-unit integrated plants, is the very large number of process simulations required. This stems from the large ensemble of process simulations needed to capture the effect of model uncertainty in conjunction with maintaining a large number of live points. We have implemented an acceleration strategy as part of nested sampling, which reduces the overall number of process model runs without impairing a design space's accuracy. A recommended follow-up to this work would entail a vectorized implementation of the nested sampling algorithm in order to further reduce the time needed to characterize a complex DS, for instance by computing simulation ensembles or replacement proposals in parallel.

Lastly, the characterization of the probabilistic DS in our work did not include adjustable control actions, although these may increase the size of the probabilistic design space. This practice is acceptable insofar as many pharmaceutical processes still use minimal online measurement and control of the CQAs---the process is instead carried to completion, followed by testing of the final product. As part of future work, it would be interesting to extend the nested sampling approach in order to encompass such recourse actions.

\begin{acknowledgement}
This work is supported by Eli Lilly and Company through the Pharmaceutical Systems Engineering Lab (PharmaSEL) program. Lucian Gomoescu is a Marie Sk\l{}odowska-Curie early stage researcher at Process Systems Entreprise Ltd enrolled in the European Union's Horizon 2020 research and innovation program under grant agreement 675585 (Marie Sk\l{}odowska-Curie ITN SyMBioSys). Radoslav Paulen gratefully acknowledges the contribution of the European Commission under grant 790017 (GuEst), the contribution of the Slovak Research and Development Agency (project APVV 15-0007), and the support of Slovak Ministry of Education, Science, Research and Sport under the project STU as the Leader of Digital Coalition 002STU-2-1/2018.. The authors are grateful to Carla Vanesa Luciani for her assistance with the Suzuki reaction case study.
\end{acknowledgement}

\paragraph{Data Statement:}
No new data was collected in the course of this research.

\begin{suppinfo}

The following supplementary information is available free of charge via the Internet at \url{https://pubs.acs.org/doi/10.1021/acs.iecr.9b05006?goto=supporting-info}.
\begin{itemize}
  \item Appendix A: Nested Sampling for Bayesian Evidence
  \item Appendix B: Michael Addition Model
  \item Appendix C: Suzuki Coupling Model
  \item Appendix D: {\sf DEUS} Python Package
\end{itemize}

\end{suppinfo}

\bibliography{references_revised}


\begin{table*}[t]
\begin{framed}
\small
\printnomenclature
\end{framed}
\end{table*}
\pagestyle{empty}

\end{document}